\documentclass[a5paper, 10pt]{article}
\UseRawInputEncoding
\usepackage{cite, amsmath, amssymb, booktabs, lastpage, graphicx}
\usepackage[hidelinks]{hyperref}
\usepackage{xcolor}
\usepackage{geometry}
\usepackage{setspace}

\usepackage{amsthm}

\theoremstyle{plain}
\newtheorem{theorem}{Theorem}
\newtheorem{conjecture}{Conjecture}

\newtheorem{problem}[theorem]{Problem}

\theoremstyle{remark}

\theoremstyle{definition}

\usepackage{graphicx}
\makeatletter
\renewcommand{\maketitle}{
	\begin{center}

		{\Large\bfseries \@title} \par
		\vspace{5mm}
		\baselineskip=0.2in
		{\large\bfseries \@author}\par
		\vspace{1mm}
		{\it \@address} \par
		{\small\tt \@email} \par
		\vspace{3mm}
		{\small (Received \@date)} \par
	\end{center}
	\vspace{3mm}
}
\newcommand{\address}[1]{\def\@address{#1}}
\newcommand{\email}[1]{\def\@email{#1}}

\makeatother

\newcommand{\acknowledgment}[1]{\vspace{5mm}\singlespacing
	{\noindent\textbf{\textit{Acknowledgement\/}:} #1}
}

\usepackage{fancyhdr}

\usepackage[margin=1cm,%
			font=footnotesize,%
			format=hang,%
			labelsep=period,%
			labelfont=bf]{caption}

\newcommand{\irr}{\,{\rm irr}}

\title{Bounds and Optimal Results\\ for the Total Irregularity Measure}
\author{Akbar Ali$^{a,}$\footnote{Corresponding author}, Darko Dimitrov$^{b}$, Tam\'as R\'eti$^{c}$,\\ Abeer M. Albalahi$^{a}$, Amjad E. Hamza$^{a}$}

\address{$^a$Department of Mathematics, College of Science,\\  University of Ha\!'il, Ha\!'il, Saudi Arabia\\
	$^b$Faculty of Information Studies, 8000 Novo Mesto, Slovenia\\
	$^c$\'Obuda University, B\'ecsi\'ut, 96/B, H-1034 Budapest, Hungary\\
    }

\email{akbarali.maths@gmail.com, darko.dimitrov11@gmail.com, reti.tamas@bgk.uni-obuda.hu, a.albalahi@uoh.edu.sa,  boaljod2@hotmail.com}

\date{XXXX}

\begin{document}

\maketitle

\begin{abstract}
A (molecular) graph in which all vertices have the same degree is known as a regular graph. According to Gutman, Hansen, and M\'elot [{\it J. Chem. Inf. Model.} {\bf45} (2005) 222--230],  it is of interest to measure the irregularity of nonregular molecular graphs both for descriptive purposes and for QSAR/QSPR studies. The graph invariants that can be used to measure the irregularity of graphs are referred to as irregularity measures. One of the well-studied irregularity measures is the ``total irregularity'' measure, which was introduced about a decade ago. Bounds and optimization problems for this measure have already been extensively studied. A considerable number of existing results (concerning this measure) also hold for molecular graphs; particularly, the ones regarding lower bounds and minimum values of the mentioned measure. The primary objective of the present review article is to collect the existing bounds and optimal results concerning the total irregularity measure. Several open problems related to the existing results on the total irregularity measure are also given.
\end{abstract}

\onehalfspacing

\section{Introduction}

All graphs considered in the present article are simple unless otherwise stated (for instance, the graphs considered in Theorem \ref{Domicolo-19-PEIS-UB-1} may not be simple). The graph-theoretical (and chemical-graph-theoretical) terminology used in this article, but not defined here, can be found in some relevant books, like \cite{Bondy08,Chartrand-16} (and \cite{Wagner-book-18,Trinajstic-book-92}, respectively).

Irregularity in graphs typically refers to the properties of graphs that deviate from certain regular patterns; for instance, see the book \cite{Ali-Book-21}.
Generally, an irregularity measure of a graph $G$ is a numerical quantity used to assess the extent to which $G$ deviates from being regular. Particularly, a numerical quantity $IM$ associated with a connected graph $G$ is said to be an irregularity measure of $G$ if it satisfies the following condition: $IM(G)\ge0$ with equality if and only if $G$ is regular; for example, see \cite{Ali-JAMC-19}.
Irregularity measures may be significant in dealing with certain problems raised in several fields, including network analysis and chemistry; for example, see \cite{Criado-14,Estrada-10,Gutman-JCIM-05,Reti-MATCH-18,Estrada-ACS-10}.

There exist many irregularity measures in the literature; for instance, see the papers \cite{Reti-MATCH-18,Dimitrov-APH-14,Gutman-KJS-16,Ma-AMC-19,Abdo-AMC-19-,Dimitrov-AMC-23} and the survey \cite{Ananias-PO-13}. In the present article, we are concerned with the total irregularity, which is an irregularity measure introduced around a decade ago by Abdo, Brandt, and Dimitrov in \cite{Abdo-IJCA-15}. The total irregularity of a graph $G$ is denoted by $\irr_t(G)$ and is defined as
$$
\irr_t(G)=\sum_{\{x,y\}\subset V}|d_x-d_y|,
$$
where $V$ is the vertex set of $G$ and $d_x$ denotes the degree of the vertex $x\in V$. The total irregularity is an extended version of Albertson's well-known irregularity measure \cite{Albertson}, defined as
$$
\irr(G)=\sum_{xy\in E}|d_x-d_y|,
$$
where $E$ is the edge set of $G$. It needs to be noted \cite{Boaventura-RMTA-15} that
\begin{equation}\label{eq-irr+irr-}
\irr_t(G)= irr(G)+ irr\left(\overline{G}\right),
\end{equation}
where $\overline{G}$ is the complement of $G$; it is defined as a graph with the same vertex set as $G$ has, while two vertices in $\overline{G}$ are adjacent if and only if they are nonadjacent in $G$.
We also note that for any (nonregular) disconnected graph $G$, we have $\irr(G) = 0$ if and only if every component of $G$ is regular. Conversely, $\irr_t(G) = 0$ holds if and only if $G$ is regular, regardless of whether $G$ is connected.

Let $X=(x_1,x_2,\dots,x_n)$ be a sequence of nonnegative real numbers such that $x_1\ge x_2\ge \cdots \ge x_n$ and $\sum_{i=1}^n x_i\ne 0$. The Gini coefficient $\zeta$ (also known as Gini index, named after Gini \cite{Gini-1912}) for $X$ can be written (see \cite{Sen-73}) as
\begin{equation}\label{Gini-index-X}
\zeta(X)=\frac{1}{2n\sum_{i=1}^n x_i} \sum_{i=1}^n\sum_{j=1}^n |x_i-x_j|=1+\frac{1}{n}-\frac{2}{n\sum_{i=1}^n x_i} \sum_{i=1}^n ix_i\,.
\end{equation}
If $G$ is an $n$-order graph with size $m\ge1$ and degree sequence $(d_1,\dots, d_n)$, such that $d_1\ge d_2\ge\cdots\ge d_n$, then by using \eqref{Gini-index-X}
the Gini index for $G$ is defined (for example, see \cite{Eliasi-IJMC-15,Ali-CM-20,Domicolo-19-PEIS}) as follows:
\begin{equation}\label{Gini-index-G}
\zeta(G)=\frac{1}{4mn} \sum_{i=1}^n\sum_{j=1}^n |d_i-d_j|=1+\frac{1}{n}-\frac{1}{nm} \sum_{i=1}^n id_i\,.    \end{equation}
Thus, $\irr_t$ of $G$ can also be rewritten as
\begin{equation*}\label{Eq-new-01}
\irr_t(G)= 2nm \cdot \zeta(G) =2(n+1)m - 2 \sum_{i=1}^n id_i\,,
\end{equation*}
where $d_1\ge d_2\ge\cdots\ge d_n$.

The main purpose of the present review article is to collect existing bounds and optimal results on $\irr_t$. The remaining part of the paper is organized as follows. The next section provides most of the notations and definitions used in the subsequent sections. Section \ref{sec-3} consists of two subsections: the first one is devoted to the results on the greatest value of $\irr_t$ over different graph classes while the second one is about its least values. Section \ref{sec-4} also is divided into two subsections: the first one gives lower bounds on $\irr_t$ while the second one provides its upper bounds.
Several open problems related to existing results on $\irr_t$ are given in Section \ref{sec-5}.

\section{Preliminaries}\label{sec-2}

In this section, we provide most of the notations and definitions used in the subsequent sections.

By an $n$-order graph, we mean a graph of order $n$. A nontrivial graph is any $n$-order graph with $n\ge2$. By a molecular graph, we mean a connected graph of the maximum degree not more than $4$. The $n$-order star, path, and complete graphs are denoted by $S_n$, $P_n$, and $K_n$, respectively. A graph that is not regular is known as a nonregular graph. A pendent vertex is a vertex with degree 1. A vertex with a degree greater than 2 is referred to as a branching vertex. The degree set of a graph $G$ is the set of all distinct elements of the degree sequence of $G$.

A nontrivial path $S$ in a tree $T$ is said to be a segment of $T$ if the end vertices of $S$ have degrees different from 2 in $T$ and every other vertex (if exists) of $S$ have degree $2$ in $T$.

The cyclomatic number of an $n$-order connected graph with size $m$ is the number $m-n+1$. By an $n$-order unicyclic, or bicyclic, or tricyclic graph, we mean a connected $n$-order graph with the cyclomatic number 1, or 2, or 3, respectively.

The distance $d(x,y)$ between two vertices $x$ and $y$ of a graph is the length of a shortest path between them. The eccentricity of a vertex $x$ in a graph with vertex set $V$ is denoted by $\epsilon_x$ and is defined as \linebreak $\epsilon_x=\max_{u\in V}d(x,u)$.
The diameter of a graph $G$ is the greatest number among the eccentricities of all vertices of $G$.
The non-self-centrality number of a connected graph $G$ is denoted as $N(G)$ and is defined \cite{Xu-AMS-16} as
\[
N(G)= \sum_{\{x,y\}\subseteq V} |\epsilon_x-\epsilon_y|.
\]

Following \cite{Ali-RMJM-24}, we recall the definitions of a polyomino system and related concepts as follows: ``A polyomino system is a connected geometric figure obtained by concatenating congruent regular squares side to side in a plane in such a way that the figure divides
the plane into one infinite (external) region and a number of finite (internal) regions,
and all internal regions must be congruent regular squares. In a polyomino system,
two squares are said to be adjacent if they share a side. The characteristic graph of a given polyomino system consists of vertices corresponding
to squares of the system; two vertices are adjacent if and only if the corresponding
squares are adjacent. A polyomino system whose characteristic graph is the path
graph is called a polyomino chain. In a polyomino chain, a square having one (respectively two) neighboring
square(s) is called a terminal (respectively, nonterminal) square. Any polyomino chain
can be represented by a graph, in which the edges represent the sides of a polyomino chain while the vertices correspond to the points where two sides of a square meet. In
what follows, by a polyomino chain, we always mean the graph corresponding to the
polyomino chain. In a polyomino chain, a nonterminal square having a vertex of degree 2
is known as a kink.  A linear polyomino chain is the one, without kinks. A zigzag polyomino chain is the one, consisting of only kinks and/or terminal squares.''

\section{Optimal Results}\label{sec-3}

This section consists of two subsections: the first one is devoted to the results about the greatest value of $\irr_t$ over different graph classes while the second one is about its least values.

\subsection{Graphs with the maximum total irregularity}

We start this section with the following result about the maximum value of $\irr_t$ in the class of all fixed-order trees:

\begin{theorem}\label{Abdo-DMTCS-14-max-thm-1} {\rm \cite{Abdo-DMTCS-14}}
The star graph $S_n$ uniquely achieves the greatest value of $\irr_t$\, among all $n$-order trees for $n\ge4$. The mentioned greatest value is $(n-1)(n-2)$.
\end{theorem}

Theorem \ref{Abdo-DMTCS-14-max-thm-1} was extended in
\cite{Eliasi-IJMC-15}, where the graphs achieving the first five greatest values of $\irr_t$\, among all $n$-order trees were determined for $n\ge13$.

The problem of finding graphs with the greatest value of $\irr_t$\, among all fixed-order trees with a given number of segments or branching vertices was addressed in \cite{Yousaf-CSF-22}.
In \cite{Yousaf-CCO-21}, the problems of finding graphs achieving the greatest value of $\irr_t$\, over the class of all fixed-order\\
(i) molecular trees, and\\
(ii) trees with a fixed maximum degree,\\
were addressed.

Next, we write about results concerning unicyclic graphs.

\begin{theorem}\label{You-AC-14-max-thm-1} {\rm \cite{You-AC-14}}
The graph constructed by inserting one edge in the star graph $S_n$\, uniquely achieves the greatest value of $\irr_t$\, among all $n$-order unicyclic graphs for $n\ge4$. The mentioned greatest value is $n^2 - n - 6 $.
\end{theorem}

Theorem \ref{You-AC-14-max-thm-1} was extended in
\cite{Eliasi-IJMC-15} where the graphs achieving the first four greatest values of $\irr_t$\, among all $n$-order unicyclic graphs were determined for $n\ge13$. Next, we write about results concerning bicyclic graphs.

\begin{theorem}\label{You-JAM-14-max-thm-1} {\rm \cite{You-JAM-14}}
The graph constructed by inserting two adjacent edges (having a common vertex) in the star $S_n$ uniquely achieves the greatest value of $\irr_t$\, among all $n$-order bicyclic graphs for $n\ge4$. The mentioned greatest value is $n^2 + n - 16 $.
\end{theorem}

The graphs achieving the first three greatest values of $\irr_t$\, among all $n$-order bicyclic graphs, were determined in
\cite{Eliasi-IJMC-15} for $n\ge12$. Particularly, Theorems \ref{Abdo-DMTCS-14-max-thm-1}, \ref{You-AC-14-max-thm-1} and \ref{You-JAM-14-max-thm-1} were proved by Eliasi \cite{Eliasi-IJMC-15} in a unified and alternative way, using the concept of the Gini index (see Equation \eqref{Gini-index-G}).

The paper \cite{Yousaf-AEJM-21} attempted to find graphs achieving the greatest value of $\irr_t$\, among all fixed-order connected graphs with cyclomatic number $\nu$ for $1\le \nu \le 6$.

A graph of size at least 1 in which both end-vertices of every edge have different degrees is known as a totally segregated graph \cite{Jackson-CN-86}. The problem of finding graphs achieving the greatest value of $\irr_t$\, over two certain classes of totally segregated bicyclic graphs of a fixed order was addressed in \cite{Jorry-CMA-21,Jorry-CMA-22}.

\begin{figure}[!ht]
 \centering
  \includegraphics[width=0.98\textwidth]{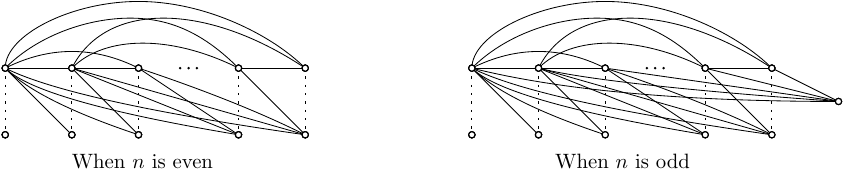}
      \caption{The $n$-order antiregular graphs $A_n$ (including the dotted edges) and $\overline{A}_n$ (excluding the dotted edges) are shown for even and odd orders, depicted on the left and right, respectively.}
    \label{Fig-1}
     \end{figure}

A nontrivial $n$-order graph
whose degree set consists of $n-1$ elements is known as an antiregular graph (or a quasi-perfect graph \cite{Bezhad-AMM-67}).
For every integer $n \ge 2$, there are exactly two nonisomorphic antiregular $n$-order graphs \cite{Bezhad-AMM-67}; one of them is connected and the other is disconnected.
Denote by $A_n$ and $\overline{A}_n$ the connected and disconnected antiregular graphs of order $n$, respectively (see Figure \ref{Fig-1}).
For $n\ge3$, let $D\subset E(A_n)$ be the set of dotted edges of $A_n$ (see Figure \ref{Fig-1}). Then $\overline{A}_n=A_n-D$. Next, we state the result concerning the maximum value of $\irr_t$ among all fixed-order graphs.

\begin{theorem}\label{Abdo-DMTCS-14-max-thm-2}{\rm \cite{Abdo-DMTCS-14}}
The graphs in \,$\{A_n-D':\emptyset\subseteq D' \subseteq D\}$ are the only graphs achieving the greatest value of $\irr_t$\, among all $n$-order graphs for every $n\ge3$. The respective greatest values are:
%
\[
\begin{aligned}
\displaystyle\frac{2 n^3 - 3 n^2 - 2 n + 3}{12} & \quad \text{when $n$ is odd, and} \\[3mm]
\displaystyle\frac{n(2 n^2 - 3n - 2)}{12} & \quad \text{when $n$ is even.}
\end{aligned}
\]

\end{theorem}

A graph is said to be a maximal $\ell$-degenerate if every subgraph has a vertex of degree at most $\ell$, and this property does not hold if any new edge is added to the graph (see \cite{Bickle-DAM-23}). The join $G_1 + G_2$ of two graphs $G_1$ and $G_2$ is the graph constructed from these two graphs by making every vertex of $G_1$ adjacent to all vertices of $G_2$. Recall that $\overline{G}$ denotes the complement of a graph $G$.

\begin{theorem}\label{Bickle-DAM-23-max-thm1} {\rm \cite{Bickle-DAM-23}}
Among all maximal $\ell$-degenerate graphs of order  $n$, the graph $K_\ell+\overline{K}_{n-\ell}$ uniquely achieves the greatest value of $\irr_t$\, for $n\ge \ell+2$. The mentioned greatest value is $\ell(n-\ell)(n-\ell-1)$.
\end{theorem}

Now, we state a result about fixed-order graphs with a given cyclomatic number and a fixed number of pendent vertices.

\begin{theorem}\label{Ghalavand-JAMC-20-UB-thm-1} {\rm \cite{Ghalavand-JAMC-20}}
In the class of all nontrivial connected $n$-order graphs with cyclomatic number $\nu$ and having $p$ pendent vertices, only the graph(s) with the degree sequence
$$
(2\nu+p, \underbrace{2,2, \ldots, 2}_{n-p-1}, \underbrace{1,1, \ldots, 1}_p\,)
$$
attain(s) the maximum value of $\irr_t$, where $n\ge p+1$. The mentioned maximum value is $(2 n-p-1)p+2(\nu-1)(n-1)$.
\end{theorem}

We end this section with the following optimal result concerning polyomino chains:

\begin{theorem}\label{Ali-16-max-thm1} {\rm \cite{Ali-JCTN-15,Yarahmadi-JNRM-21}}
   Among all polyomino chains consisting of $s$ squares, the zigzag polyomino chain uniquely achieves the greatest value of $\irr_t$\, for every $s\ge3$. The mentioned greatest value is $2(s^2 + 2s - 4)$.
\end{theorem}

\subsection{Graphs with the minimum total irregularity}

Most of the results of this section also hold for molecular graphs; particularly, it is true when the graph(s) attaining the minimum value of $\irr_t$ in the considered class of graphs is/are molecular.
We start writing about trees having the minimum value of $\irr_t$.

\begin{theorem}\label{Zhu-Filomat-16-min-thm-1} {\rm \cite{Zhu-Filomat-16}}
The path graph $P_n$ uniquely achieves the least value of $\irr_t$\, among all $n$-order trees for every $n\ge4$. The mentioned least value is $2(n-2)$.
\end{theorem}

In \cite{Zhu-Filomat-16}, the graphs achieving the second- and third-least values of $\irr_t$ among all $n$-order trees were also determined for $n\ge6$. Ghalavand and Ashrafi \cite{Ghalavand-JAMC-20} extended these results by determining the fourth to seventh least values of $\irr_t$ over the class of all fixed-order trees.
Moreover, Theorem \ref{Zhu-Filomat-16-min-thm-1} was proved by Eliasi \cite{Eliasi-IJMC-15} (see also \cite{Ghalavand-JAMC-20}) in an alternative way, using the concept of the Gini index, see Equation \eqref{Gini-index-G}.

The problem of finding graphs with the least value of $\irr_t$\, among all fixed-order trees with a given number of segments or branching vertices was addressed in \cite{Yousaf-CSF-22}.

Recall that a graph of size at least 1 in which both end-vertices of every edge have different degrees is known as a totally segregated graph \cite{Jackson-CN-86}.
The problem of finding graphs achieving the least value of $\irr_t$\, among all totally segregated trees of a fixed order was attacked in \cite{Jorry-IJMTT-15}.

Next, we write about the results concerning unicyclic graphs.
The problem of determining the graph attaining the minimum value of $\irr_t$ among all fixed-order unicyclic graphs has a trivial solution (that is, the cycle graph); the next result gives a solution to a similar problem when one considers all nonregular unicyclic graphs.

\begin{theorem}\label{Zhu-Filomat-16-min-thm-2} {\rm \cite{Zhu-Filomat-16}}
Only the graphs with the degree sequence $(3,2,\dots,2,1)$ achieve the least value of $\irr_t$\, among all $n$-order nonregular unicyclic graphs for every $n\ge5$. The mentioned least value is $2(n-1)$.
\end{theorem}

In \cite{Zhu-Filomat-16}, the graphs achieving the second-minimum value of $\irr_t$ among all $n$-order nonregular unicyclic graphs were also found for every $n\ge5$.
Ghalavand and Ashrafi \cite{Ghalavand-JAMC-20} extended these results by determining the graphs achieving the third-minimum value of $\irr_t$ among all $n$-order nonregular unicyclic graphs for every $n\ge7$.

Now, we write regarding the results on bicyclic graphs.

\begin{theorem}\label{Zhu-Filomat-16-min-thm-3} {\rm \cite{Zhu-Filomat-16}}
Only the graphs with the degree sequence $(3,3,\overbrace{2,\dots,2}^{n-2}\,)$ achieve the least value of $\irr_t$\, among all $n$-order bicyclic graphs for every $n\ge7$. The mentioned least value is $2(n-2)$.
\end{theorem}

The graphs achieving the second- and third-minimum values of $\irr_t$ among all $n$-order bicyclic graphs were also found in \cite{Zhu-Filomat-16} for every $n\ge7$.
Ghalavand and Ashrafi \cite{Ghalavand-JAMC-20} extended these results by determining the graphs achieving the fourth-minimum value of $\irr_t$ among all $n$-order bicyclic graphs for every $n\ge9$.

The paper \cite{Jorry-IJIT-19} addressed the problem of finding graphs achieving the least value of $\irr_t$\, over a certain class of totally segregated bicyclic graphs of a fixed order.

Next, we write about the results regarding tricyclic graphs.

\begin{theorem}\label{Ahmed-KJS-23-min-thm-1} {\rm \cite{Ahmed-KJS-23}}
Only the graphs with the degree sequence $$(3,3,3,3,\overbrace{2,\dots,2}^{n-4}\,)$$ achieve the least value of $\irr_t$\, among all $n$-order tricyclic graphs for every $n\ge7$. The mentioned least value is $4(n-4)$.
\end{theorem}

The problem of determining graphs achieving the second- and third-minimum values of $\irr_t$ among all $n$-order tricyclic graphs was also attacked in \cite{Ahmed-KJS-23} for $n\ge7$; however, the degree sequences mentioned in both cases are not graphical.

Next, we state two results about fixed-order graphs with a given cyclomatic number.

\begin{theorem}\label{Ghalavand-JAMC-20-Min-thm-1} {\rm \cite{Ghalavand-JAMC-20}}
  Among all $n$-order molecular graphs with cyclomatic number $\nu$, the graphs with degree sequence
$$
(\,\underbrace{4,4, \ldots, 4}_{k-1}, \underbrace{3,3, \ldots, 3}_{2(\nu+k-2)} ,\underbrace{2,2, \ldots, 2}_{n-2\nu-3k+5})
$$
achieve the $k$-th least value of $\irr_t$\,, where $n>2 \nu^2-3\nu+2$, $\nu\geq 2$ and $1 \leq k \leq \nu$.
\end{theorem}

\begin{theorem}\label{Ghalavand-JAMC-20-Min-thm-2} {\rm \cite{Ghalavand-JAMC-20}}
  Among all $n$-order connected graphs with cyclomatic number $\nu$, the graphs with degree sequence
$$
(\,\underbrace{4,4, \ldots, 4}_{k-1}, \underbrace{3,3, \ldots, 3}_{2(\nu+k-2)} ,\underbrace{2,2, \ldots, 2}_{n-2\nu-3k+5})
$$
achieve the $k$-th least value of $\irr_t$\,, where $n>2 \nu^2-3\nu+2$, $\nu\geq 3$ and $1 \leq k \leq 3$.
\end{theorem}

Now, we state two results about fixed-order graphs with a given cyclomatic number and a fixed number of pendent vertices.

\begin{theorem}\label{Ghalavand-JAMC-20-LB-thm-1} {\rm \cite{Ghalavand-JAMC-20}}
In the class of all nontrivial connected $n$-order graphs with cyclomatic number $\nu$ and having $p$ pendent vertices, only the graph(s) with the degree sequence
$$
(\,\underbrace{3,3, \ldots, 3}_{p+2(\nu-1)}, \underbrace{2,2, \ldots, 2}_{n-2(p+\nu-1)}, \underbrace{1,1, \ldots, 1}_p\,)
$$
attain(s) the minimum value of $\irr_t$, where $n\ge2(p+\nu-1)$. The mentioned minimum value is
\[
\frac{1}{2}(n+2\nu-2)(n-2\nu+2)-2\left(p+\nu-\frac{1}{2} n-1\right)^2.
\]
\end{theorem}

\begin{theorem}\label{Ghalavand-JAMC-20-LB-thm-2} {\rm \cite{Ghalavand-JAMC-20}}
In the class of all connected $n$-order graphs of maximum degree at least $5$, with cyclomatic number $\nu$ and having $p$ pendent vertices, only the graph(s) with the degree sequence
$$
(\,5,\underbrace{3, 3,\ldots, 3}_{p+2\nu-5}, \underbrace{2,2, \ldots, 2}_{n-2(p+\nu-2)} \underbrace{1,1, \ldots, 1}_p\,)
$$
attain(s) the minimum value of $\irr_t$, where $n\ge2(p+\nu-2)$. The mentioned minimum value is
\[
2[ \nu(n-2 \nu+8)+ p(n-2\nu-p+4)-(n+9)].
\]
\end{theorem}

Next, we turn our attention toward the minimum value of $\irr_t$\, over the class of connected fixed-order nonregular graphs.

\begin{conjecture}\label{Zhu-Filomat-16-min-conj} {\rm \cite{Zhu-Filomat-16}}
   The least value of $\irr_t$\, among all connected $n$-order nonregular graphs (with $n\ge3$) is
   \[
\begin{aligned}
n-1,  & \quad \text{when $n$ is odd, and}\\[2mm]
2(n-2)  & \quad \text{when $n$ is even.}
\end{aligned}
\]
\end{conjecture}
Abdo and Dimitrov  \cite{Abdo-BAMS-15} proved Conjecture \ref{Zhu-Filomat-16-min-conj} and determined the second-minimum value of $\irr_t$\, among all connected $n$-order nonregular graphs. Ashrafi and Ghalavand \cite{Ashrafi-AMC-20} proved Conjecture \ref{Zhu-Filomat-16-min-conj} in an alternative but a shorter way, and determined further ordering concerning the least value of $\irr_t$\, over the class of all connected $n$-order nonregular graphs.



An Eulerian graph is a connected graph of order at least 3 in which every vertex has an even degree. Abdo and Dimitrov  \cite{Abdo-BAMS-15} proved that the inequality $\irr_t(G)\ge n-1$ holds for every connected $n$-order nonregular graph $G$, where the equality $\irr_t(G)= n-1$ is achieved only when $n$ is odd, and the graphs attaining this equality are non-Eulerian. Thus, by \cite{Abdo-BAMS-15}, the inequality
\begin{equation}\label{eq-AA-00001}
\irr_t(G)\ge 2(n-2)
\end{equation}
holds for every Eulerian nonregular connected graph $G$ of order $n\ge5$. Using the bound \eqref{eq-AA-00001}, Nasiri et al. \cite{Nasiri-IJMC-18} obtained the next two results.

\begin{theorem}\label{Nasiri-IJMC-18-min-thm1} {\rm \cite{Nasiri-IJMC-18}}
   Among all those connected graphs of order at least $5$ that are nonregular and Eulerian, the graph $S^\star_5$\, constructed by adding two nonadjacent edges to the $5$-order star graph $S_5$ uniquely achieves the least value of $\irr_t$. The mentioned least value is $8$.
\end{theorem}

\begin{figure}[!ht]
 \centering
\includegraphics[width=0.5\textwidth]{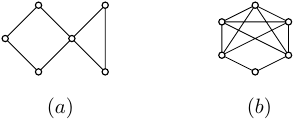}
   \caption{Two graphs referred in Theorem \ref{Nasiri-IJMC-18-min-thm2}.}
    \label{Fig-2}
     \end{figure}

\begin{theorem}\label{Nasiri-IJMC-18-min-thm2} {\rm \cite{Nasiri-IJMC-18}}
   Among all those connected graphs of order at least $5$ that are nonregular, Eulerian, and different from $S^\star_5$ (defined in Theorem $\ref{Nasiri-IJMC-18-min-thm1}$), the graphs shown in Figure $\ref{Fig-2}$ are the only graphs achieving the least value of $\irr_t$. The mentioned least value is $10$.
\end{theorem}

We close this section by stating the following optimal result about polyomino chains:

\begin{theorem}\label{Ali-16-min-thm1} {\rm \cite{Yarahmadi-JNRM-21,Ali-JCTN-15}}
   Among all polyomino chains consisting of $s$ squares, the linear polyomino chain uniquely achieves the least value of $\irr_t$\, for every $s\ge3$. The mentioned least value is $8(s-1)$.
\end{theorem}

\section{Bounds on the Total Irregularity}\label{sec-4}
This section is divided into two subsections: the first one gives lower bounds on $\irr_t$ while the second one provides its upper bounds. Before going towards the first subsection of the present section, we remark here that we can derive in many cases bounds concerning $\irr_t$ from the existing bounds of $\irr$; for instance, for a graph $G$ of size $m$, the chain of inequalities
\[
\sqrt{\sigma(G)}  \le \irr(G) \le \sqrt{m\cdot\sigma(G)},
\]
holds (see \cite{Matejic-SPSUNP-19}) with either of the equalities if and only if $G$ is regular, where $\sigma(G)$ is the sigma index \cite{Furtula-BASA-15} of $G$ (see also \cite{Ali-DAM-23,Gutman-MATCH-18,Dimitrov-DML-23,kpvsd-sict-2024}) defined as
\[
\sigma(G) =\sum_{xy\in E} (d_x-d_y)^2;
\]
consequently, for an $n$-order graph $G$ of size $m$,
by Equation \eqref{eq-irr+irr-}, we have
\[
\sqrt{\sigma(G)}+\sqrt{\sigma\left(\overline{G}\right)}  \le \irr_t(G) \le \sqrt{m\cdot\sigma(G)}+\sqrt{\left(\frac{n(n-1)}{2}-m\right)\sigma\left(\overline{G}\right)},
\]
with either of the equalities if and only if $G$ is regular.

\subsection{Lower bounds on the total irregularity}

We start this section with the remark that the inequality $\irr_t(G)\ge \irr(G)$ holds for every graph $G$, because of the definitions of $\irr_t$ and $\irr$.

Let $n_s$ denote the number of vertices of a graph $G$ having degree $s$.

\begin{theorem}\label{Ashrafi-AMC-20-LB-1}{\rm \cite{Ashrafi-AMC-20}}
If $G$ is an $n$-order nonregular graph, then $$\irr_t(G) \ge n_s(n -n_s),$$ where $n_s\ne 0$ with $1 \le s \le n - 1$.
\end{theorem}

The first Zagreb index\, ---\, introduced in \cite{Gutman-72} (see its most-recent surveys \cite{Borovicanin-17-MATCH,Ali-MATCH-18})\, --- \,for a graph $G$ is denoted by $Z_1(G)$ and is defined as
\[
Z_1(G)= \sum_{x\in V} (d_x)^2.
\]

\begin{theorem}\label{Ghalavand-23-JM-LB-thm-1} {\rm \cite{Ghalavand-23-JM}}
Let $G$ be a nonregular $n$-order graph with maximum degree $\Delta$, size $m$ and minimum degree $\delta$. Then,
$$
\irr_t(G) \ge \frac{n\cdot Z_1(G)-4 m^2}{\Delta-\delta},
$$
with equality if and only if the degree set of $G$ is $\{\delta, \Delta\}$.
Also,
$$
\irr_t(G) \geq \frac{n^2 \Delta(\Delta-1)+2 m n(2 \Delta-1)-4 m^2}{\Delta-\delta},
$$
with the equality if and only if the degree set of $G$ is $\{\Delta-1,\Delta\}$. Furthermore,
$$
\irr_t(G) \geq \frac{2 m n(2 \delta+1)-4 m^2-n^2 \delta(\delta+1)}{\Delta-\delta},
$$
with equality if and only if the degree set of $G$ is $\{\delta, \delta+1\}$.
\end{theorem}

A graph in which all vertices have the same eccentricity is known as a self-centered graph. By a non-self-centered graph, we mean a graph that is not self-centered.

Let $K_n(k)$ denote the $n$-order graph constructed from the complete graph $K_n$ by dropping $k$ pairwise nonadjacent edge(s), where $n \ge 4$ and
$1 \le k \le 	\lfloor n/2 \rfloor$. Take $\mathcal{G}_n = \{K_n(k) : 1 \le k \le 	\lfloor n/2 \rfloor
\}$.

\begin{theorem}\label{Xu-AMC-18-LB-1} {\rm \cite{Xu-AMC-18}}
If $G$ is an $n$-order connected non-self-centered graph with diameter $2$, then $\irr_t(G) \ge N(G)$ with equality if and only if $G \in \mathcal{G}_n$, where $n\ge3$ and $N(G)$ is the non-self-centrality number of $G$ defined in Section $\ref{sec-2}$.
\end{theorem}

\begin{theorem}\label{Xu-AMC-18-LB-2} {\rm \cite{Xu-AMC-18}}
Let $G \not\in \mathcal{G}_n$ be an $n$-order connected non-self-centered graph with diameter $2$, where $n\ge3$. If $n'$ is the number
of those vertices of $G$ that have degree less than $n - 1$, then
\[
\irr_t(G)>
\begin{cases}
N(G) + 4(n-3)  & \text{when $n'$ is even,}\\[2mm]
N(G) + 2(n-2)  & \text{when $n'$ is odd.}
\end{cases}
\]
\end{theorem}

\begin{theorem}\label{Xu-AMC-18-LB-3} {\rm \cite{Xu-AMC-18}}
If $G$ is a connected nonregular graph with a diameter of at least $4$, then either
$\irr_t(G)> N(G)$ or $\irr_t(G)> N\left(\overline{G}\right)$.
\end{theorem}

In Theorem \ref{Xu-AMC-18-LB-3}, if the text ``diameter of at least 4'' is replaced with ``diameter of at least 3'' then the resulting statement does not hold generally; however, the following holds:

\begin{theorem}\label{Xu-AMC-18-LB-4} {\rm \cite{Xu-AMC-18}}
For every integer $d$ greater than $2$, there exists at least one connected graph $G$ with order at least $6$ and diameter $d$ such that $\irr_t(G) > N (G)$.
\end{theorem}

The next result improves Theorem \ref{Xu-AMC-18-LB-4}.

\begin{theorem}\label{Tang-AMC-19-LB-1} {\rm \cite{Tang-AMC-19}}
For every integer $d$ greater than $1$, there exists at least one connected non-self-centered graph $G$ with order at least $6$ and diameter $d$ such that $\irr_t(G) > N (G)$.
\end{theorem}

\begin{problem}\label{Xu-AMC-18-LB-prblm} {\rm \cite{Xu-AMC-18}}
    Determine the class $\mathcal{G}$ of all connected graphs of diameter at least $3$ such that $\irr_t(G)>N(G)$ for every $G\in\mathcal{G}$.
\end{problem}

The next two results provide a partial solution to Problem \ref{Xu-AMC-18-LB-prblm}.

\begin{theorem}\label{Tang-AMC-19-LB-2} {\rm \cite{Tang-AMC-19}}
Let $G$ be a connected graph of diameter $3$. If $G$ is a tree different from the $4$-order path graph $P_4$ or if $G$ is a unicyclic graph different from the $6$-order cycle graph $C_6$, then $\irr_t(G) > N (G)$.
\end{theorem}

A vertex $v$ in a connected graph $G$ is said to be a central vertex of $G$ if the eccentricity of $v$ is the same as the radius of $G$.

\begin{theorem}\label{Tang-AMC-19-LB-3} {\rm \cite{Tang-AMC-19}}
Let $T$ be a tree of diameter $4$ and order at least $8$. If $u$ is a central vertex of $T$ such that $u$ has $\alpha$ neighbors of degree at least $2$ and $3\le \alpha < d_u\,$, then $\irr_t(G) > N (G)$.
\end{theorem}

\begin{theorem}\label{Ghalavand-DML-20-LB-1} {\rm \cite{Ghalavand-DML-20}}
If $G$ is a nontrivial $n$-order graph, then
$$
\irr_t(G)\ge \frac{1}{2(n-1)}\left|\operatorname{irr}_{t, 2}(G)-\operatorname{irr}_{t, 2}(\overline{G})\right| ,
$$
where $$
\operatorname{irr}_{t, 2}=\sum_{\{u, v\} \subseteq V}\left|(d_u)^2-(d_v)^2\right|.
$$
Also, if $G$ is a connected nonregular $n$-order with the condition $n \geq 3$, then
\begin{equation}\label{eq-Ghal-DML}
\irr_t(G)\ge \frac{1}{2n-3}\,\irr_{t, 2}(G),
\end{equation}
where the equality in \eqref{eq-Ghal-DML} holds when the degree set of $G$ is $\{n-2, n-1\}$.
\end{theorem}

By a universal vertex in a graph $G$, we mean a vertex adjacent to all other vertices of $G$.

\begin{theorem}\label{Lin-CM-21-LB-1} {\rm \cite{Lin-CM-21}}
If $G$ is a nonregular connected graph with radius $r$, then
$
 \irr_t(G) \ge r\cdot HA(G),
$
where $HA(G)$ is the Harary-Albertson index of $G$ defined as
$$
HA(G)=\sum_{\{u, v\} \subseteq V} \frac{|d_u-d_v|}{d(u, v)}.
$$
In addition, if $G$ does not have any universal vertex, then
$$
\irr_t(G) \ge HA(\overline{G}) +\frac{1}{2} \irr(G) .
$$
\end{theorem}

For the next result, recall the definition of a maximal $\ell$-degenerate graph given right before Theorem $\ref{Bickle-DAM-23-max-thm1}$.

\begin{theorem}\label{Bickle-DAM-LB-thm2} {\rm \cite{Bickle-DAM-23}}
Let $G$ be a maximal $2$-degenerate graph.
\begin{description}
  \item[(i).] ~\,If the order of $G$ is more than $5$ and if it contains at least two pendent vertices, then $\irr_t(G) \ge 10$.
  \item[(ii).] \,If the maximum degree of $G$ is at least $5$, then $\irr_t(G) \ge 9$.
  \item[(iii).] If $G$ has at least $7$ vertices, then $\irr_t(G) \ge 8$.
\end{description}
\end{theorem}

We end this section with the remark that lower bounds on $\irr_t$ for graphs under several graph operations can be found in
\cite{Dehgardi-AMN-20,Azari-CCO-24}.

\subsection{Upper bounds on the total irregularity}

In this section, we state the existing upper bounds on $\irr_t$.

\begin{theorem}\label{Dimitrov-AMC-15-UB-1} {\rm \cite{Dimitrov-AMC-15}}
For every connected $n$-order graph $G$, the following inequality holds:
\begin{equation}\label{Dimitrov-AMC-15-UB-1-eq}
\irr_t(G)\le \frac{n^2}{4}\cdot irr(G).
\end{equation}
The equality in \eqref{Dimitrov-AMC-15-UB-1-eq} is attained for infinitely many values of $n$.
\end{theorem}

The inequality given in Theorem \ref{Dimitrov-AMC-15-UB-1} can be improved if we restrict ourselves to trees.

\begin{theorem}\label{Dimitrov-AMC-15-UB-2} {\rm \cite{Dimitrov-AMC-15}}
For a nontrivial $n$-order tree $T$, the following inequality holds:
\begin{equation*}\label{Dimitrov-AMC-15-UB-2-eq}
\irr_t(T)\le (n-2)\, irr(T),
\end{equation*}
with equality if and only if $T$ is the path graph $P_n$.
\end{theorem}

\begin{theorem}\label{Domicolo-19-PEIS-UB-1} {\rm \cite{Domicolo-19-PEIS}}
For an $n$-order graph $G$, possibly containing loops, with $m$ edges, the following inequality holds:
\begin{equation}\label{Domicolo-19-PEIS-eq1}
\irr_t(G)\le 2(n-1)m.
\end{equation}
If $G$ contains $n-1$ vertices of degree $0$ and one vertex incident with $m$ loops then the equality in \eqref{Domicolo-19-PEIS-eq1} is attained.
\end{theorem}

The Cartesian product $G_1 \square G_2$ of two graphs $G_1$ and $G_2$ is the graph with the vertex set $V (G_1) \times V (G_2)$, where two vertices $(a_1,a_2)$ and $(b_1,b_2)$ of $G_1 \square G_2$ are adjacent if and only if either $a_1b_1$ is an edge in $G_1$ and $a_2=b_2$, or $a_1 = b_1$  and $a_2b_2$ is an edge in $G_2$.

\begin{theorem}\label{Xu-AMC-18-UB-1} {\rm \cite{Xu-AMC-18}}
If $T$ is an $n$-order tree with diameter $d$  and maximum degree $3$, then
$$\irr_t(T)<N(T)\quad \text{and} \quad \irr_t(T\square K_2)<N(T\square K_2),$$
where $n \ge15$, $N(G)$ is the non-self-centrality number of $G$ defined in Section $\ref{sec-2}$ and $d \ge\frac{2}{3}n$.
\end{theorem}

For a graph $G$, let $G^*$ be the graph constructed from $G$ by attaching one pendent vertex at every vertex of $G$.

\begin{theorem}\label{Xu-AMC-18-UB-2} {\rm \cite{Xu-AMC-18}}
If $n\ge14$ then
$$\irr_t(P_n^*)<N(P_n^*)\quad \text{and} \quad \irr_t((P_n^*)^*)<N((P_n^*)^*),$$
where $P_n$ is the $n$-order path graph.
\end{theorem}

\begin{problem}\label{Xu-AMC-18-UB-prblm} {\rm \cite{Xu-AMC-18}}
    Determine the class $\mathcal{T}$ of all trees of maximum degree at most $4$ such that the inequality $\irr_t(T)<N(T)$ holds for every $T\in\mathcal{T}$.
\end{problem}

The next result gives a partial solution to Problem \ref{Xu-AMC-18-UB-prblm}.

\begin{theorem}\label{Tang-AMC-19-UB-1} {\rm \cite{Tang-AMC-19}}
If $T$ is an $n$-order tree with diameter $d$  and maximum degree $4$ such that $G$ does not contain any vertex of degree $3$, then
$$\irr_t(T)<N(T),$$ where $n \ge10$ and $d \ge\frac{\sqrt{26}+2}{11}n$.
\end{theorem}

\begin{theorem}\label{Ghalavand-DML-20-UB-1} {\rm \cite{Ghalavand-DML-20}}
If $G$ is a connected nonregular graph of order at least $3$, then
\begin{equation}\label{eq-Ghal-DML-01}
\irr_t(G)\le \frac{1}{3}\,\irr_{t, 2}(G),
\end{equation}
where $\irr_{t, 2}$ is defined in Theorem $\ref{Ghalavand-DML-20-LB-1}$. The equality in \eqref{eq-Ghal-DML-01} holds when $G$ is a path graph.
\end{theorem}

\begin{theorem}\label{Lin-CM-21-UB-1} {\rm \cite{Lin-CM-21}}
If $G$ is a nonregular connected graph with diameter $d$, then
$ \irr_t(G) \le d\cdot HA(G),$
where $HA(G)$ is the Harary-Albertson index of $G$ defined in Theorem $\ref{Lin-CM-21-LB-1}$.
Also,
$$ \irr_t(G) \le
\begin{cases}
3\cdot HA(\overline{G})    &  \text{if $d\ge3$,}\\[2mm]
2\cdot HA(\overline{G}) &  \text{if $d\ge4$.}
\end{cases}
$$

\end{theorem}

Upper bounds on $\irr_t$ for graphs concerning several graph operations can be found in
\cite{Abdo-MMN-14}. Also, upper bounds
on $\irr_t$ for several types of composite graphs can be found in \cite{Abdo-IJCA-15}.

\section{Open Problems}\label{sec-5}

Although numerous open problems concerning $\irr_t$ can be posed, the problems presented here are closely related to existing results on $\irr_t$.

Open problems and conjectures concerning the closely related irregularity measures, $\sigma$-irregularity and $\sigma_t$-irregularity, will be omitted here. For some of these problems and conjectures, we refer the reader to \cite{Dimitrov-DML-23, Dimitrov-AMC-23, fdks-srosti-2024, ksfd-eagmtsi-2024}.

The problem of finding graphs with the extremum values of $\irr_t$\, among all fixed-order trees with a given number of segments or branching vertices was addressed in \cite{Yousaf-CSF-22}; in the case of the maximum value of $\irr_t$, the obtained graphs are not molecular ones. This observation suggests the following problem:

\begin{problem}
Characterize graphs attaining the maximum value of $\irr_t$\, in the class all fixed-order molecular trees with a given number of (i) segments or (ii) branching vertices.
\end{problem}

Theorem \ref{Abdo-DMTCS-14-max-thm-2} provides all graphs possessing the maximum value of $\irr_t$ over the class of all fixed-order connected graphs; one may think of establish a similar result by fixing the size instead of order, which yields the next problem.

\begin{problem}\label{prb-size}
    Characterize graphs attaining the maximum value of $\irr_t$ over each of the following graph classes: (i) class of all connected graphs with a given size, (ii) class of all molecular graphs with a given size.
\end{problem}

It is natural to ask establishing a maximal version of Theorem \ref{Ghalavand-JAMC-20-Min-thm-1}; or, at least establishing the mentioned version about the first three maximum values of $\irr_t$. Thus, we pose the following problem related Theorem \ref{Ghalavand-JAMC-20-Min-thm-1}:

\begin{problem}
Characterize graphs achieving the first three maximum values of $\irr_t$  among all fixed-order molecular graphs with a given cyclomatic number.
\end{problem}

Keeping in mind Theorem \ref{Ghalavand-JAMC-20-Min-thm-2} and existing literature, we pose the next problem.

\begin{problem}
Characterize graphs achieving the first three maximum values of $\irr_t$  among all fixed-order connected graphs with a given cyclomatic number (greater than $6$).
\end{problem}

Although Theorems \ref{Tang-AMC-19-LB-2} and \ref{Tang-AMC-19-LB-3} provide a partial solution to Problem \ref{Xu-AMC-18-LB-prblm}, and Theorem \ref{Tang-AMC-19-UB-1} gives a partial solution to Problem \ref{Xu-AMC-18-UB-prblm}, both of these problems are still generally open.

We end this paper by posing the minimal version of Problem \ref{prb-size} as follows:

\begin{problem}
    Characterize graphs attaining the minimum value of $\irr_t$ over each of the graph classes given in Problem $\ref{prb-size}$.
\end{problem}

\acknowledgment{This research has been funded by the Scientific Research Deanship at the University of Ha\!'il - Saudi Arabia through
project number RG-24\,046.}

\singlespacing

\end{document}